\documentclass[12pt]{amsart}%
\usepackage{amsmath,amsthm,amsfonts,amssymb}
\usepackage{verbatim}
\usepackage{hyperref}
\usepackage{latexsym}
\usepackage{marginnote}
\usepackage{amscd}
\usepackage{xcolor}
\usepackage{multicol}

\usepackage{mathrsfs}
\usepackage{graphicx}%
\providecommand{\U}[1]{\protect\rule{.1in}{.1in}}
\providecommand{\U}[1]{\protect\rule{.1in}{.1in}}
\providecommand{\U}[1]{\protect\rule{.1in}{.1in}}
\newtheorem{theorem}{Theorem}[section]
\newtheorem{lemma}[theorem]{Lemma}
\newtheorem{proposition}[theorem]{Proposition}
\newtheorem{definition}[theorem]{Definition}

\newtheorem{example}[theorem]{Example}

\begin{document}
\title[Sobolev  spaces]{Sobolev spaces arising from the Stockwell transform on  Gelfand pairs}
\author{Claude G. Dosseh, Mawoussi Todjro and  Yaogan Mensah}
\address{Department of Mathematics, University of Lom\'e, Togo}
\email{\textcolor[rgb]{0.00,0.00,0.84}{claudeganyo@gmail.com}}
\address{Department of Mathematics, University of Kara, Togo}
\email{\textcolor[rgb]{0.00,0.00,0.84}{todjrom7@gmail.com}}
\address{Department of Mathematics, University of Lom\'e, Togo}
\email{\textcolor[rgb]{0.00,0.00,0.84}{mensahyaogan2@gmail.com}}

\begin{abstract}
We introduce a family of Sobolev-type spaces on Gelfand pairs, built from the Stockwell transform: a function belongs to the space of a given smoothness when its Stockwell transform, multiplied by a weight that grows with the smoothness, is square integrable. We prove that these spaces are Hilbert spaces, that they form a decreasing scale under continuous embeddings, and that an interpolation inequality controls the norm at an intermediate smoothness by the norms at the two endpoints. Under integrability conditions on the weight and on the family of analysing windows, they embed continuously into Lebesgue spaces. Under a local domination condition on the window, we show that their elements have bounded continuous representatives. On the Euclidean motion group, where the Stockwell transform reduces to a windowed Hankel transform, this gives boundedness and continuity as soon as the smoothness exceeds half the dimension, which is the classical threshold for Euclidean spaces.

\end{abstract}
\maketitle


\section{Introduction}
  The Stockwell transform, introduced by Stockwell, Mansinha and Lowe \cite{Stockwell},  can be viewed either as a short-time Fourier transform with a frequency-dependent window, or as a phase-corrected wavelet transform, and it has since become a standard tool for analysing non-stationary signals in seismology, biomedical imaging, and mechanical fault diagnosis.

While the Stockwell transform was first developed on the real line and on Euclidean space, subsequent work has sought to place it within the broader framework of abstract harmonic analysis, in the same way the short-time Fourier and wavelet transforms have already been extended to locally compact groups and homogeneous spaces. Esmaeelzadeh carried out the first such extension \cite{Esmaeelzadeh2023, Esmaeelzadeh2024}, defining the Stockwell transform, together with its associated localization operators, on locally compact abelian groups. More recently, this construction was pushed further, to the setting of Gelfand pairs \cite{Dosseh2026}, a class that strictly contains locally compact abelian groups and that includes many spaces of independent geometric interest, such as spheres, hyperbolic spaces, and other symmetric spaces. Gelfand pairs occupy a distinguished position in noncommutative harmonic analysis;  although the ambient group need not be commutative, the algebra of bi-invariant functions still is, so that a  spherical Fourier analysis remains available. Working at the level of Gelfand pairs therefore makes it possible to transport time-frequency methods, originally designed for commutative groups, to a much wider and geometrically richer class of settings.

On the other hand, Sobolev spaces have long played a central role in the study of partial differential equations and integral equations, providing the natural framework in which questions of regularity, boundedness, and embedding are addressed. The interplay between Sobolev-type smoothness and time-frequency transforms has proved fruitful. Classically,   Sobolev spaces of Bessel potentiel type are built from the Fourier transform on Euclidean spaces \cite{Adams2003} and on locally compact groups/hypergroups \cite{Gorka2014, Kumar} and  homogeneous spaces \cite{ Bataka2025,  Krukowski,  Mensah}. 
 
The purpose of the present paper is to  construct a family of Sobolev-type spaces associated with the Stockwell transform on Gelfand pairs, and to study some of their major   properties. We show that these spaces carry a natural Hilbert space structure, we establish a scale of continuous embeddings relating spaces of different smoothness to one another and to the underlying square-integrable functions, and we obtain an embedding into  Lebesgue spaces under a suitable integrability condition. We further show, under a  local domination hypothesis, that sufficiently smooth elements of these spaces are in fact bounded and continuous. Throughout, the theory is illustrated on the Euclidean motion group, a  noncommutative Gelfand pair on which it recovers, as a special case, the classical smoothness threshold for the  continuous embedding into bounded continuous functions on  Euclidean spaces.

The remainder of the paper is organized as follows. Section \ref{prelim} collects the background material on harmonic analysis on Gelfand pairs and on the Stockwell transform needed in the sequel. Section \ref{main:results} contains the main results.

\section{Preliminaries}\label{prelim}
In this section, we provide the tools needed in the sequel.
\subsection{Harmonic analysis on Gelfand pairs}
For this subsection we refer to \cite{Dijk2009, Wolf}.

Let $G$ be a locally compact group  (assumed to be Hausdorff once for all) with a left Haar measure and let $K$  be a compact subgroup of $G$ endowed with its normalized Haar measure. A function $f: G \to \mathbb{C}$ is said to be  $K$-bi-invariant if it satisfies the property : $\forall x\in G, \forall k, k' \in K,$ 
$$f(kxk')=f(x).$$

Let $\mathcal{C}_c(G//K)$ and $L^p(G//K),\, 1\le p<\infty,$ denote respectively the space of complex-valued continuous compactly supported $K$-bi-invariant functions on $G$ and the space of complex-valued $K$-bi-invariant Haar $p$-integrable functions on $G$.
For $1\leq p<\infty$, $L^p(G//K)$ is a complex Banach space with respect to the norm 
$$\|f\|_{L^p}=\left(\int_G|f(x)|^p\right)^{\frac{1}{p}}.$$

  The convolution of $f,g \in L^1(G//K)$ is defined by
$$f\ast g(x)=\int_G f(y)g(y^{-1}x)dy,\quad x\in G.$$

The pair $(G, K)$ is called a Gelfand pair if $\mathcal{C}_c(G//K)$ is a commutative convolution algebra; by density of $\mathcal{C}_c(G//K)$ in $L^1(G//K)$, this is equivalent to say that $L^1(G//K)$ is a commutative convolution algebra.  It is well known that if $(G,K)$ is a Gelfand pair, then $G$ is unimodular \cite[Proposition 6.1.2]{Dijk2009}. 

A continuous $K$-bi-invariant function $\varphi$ is called a spherical function for the Gelfand pair $(G,K)$ if the map $\chi_\varphi: \mathcal{C}_c(G//K)\to \mathbb{C}, f\mapsto \displaystyle\int_G f(x)\varphi(x^{-1})dx$ is a nontrivial character of the convolution  algebra $\mathcal{C}_c(G//K)$, that is 
$$\chi(f\ast g)=\chi (f)\chi (g),\quad f,g\in \mathcal{C}_c(G//K).$$

The set of bounded spherical functions for the Gelfand pair $(G, K)$ is denoted by $\mathcal{S}_b(G, K)$. 
The spherical Fourier transform for $f\in L^1(G//K)$ is defined by 
$$\widehat{f}(\varphi)=\int_G f(x)\varphi(x^{-1})dx,\quad \varphi \in \mathcal{S}_b(G, K).$$

A function $\varphi: G\to \mathbb{C}$ is said to be  positive-definite function if
$$
\sum_{n=1}^{N}\sum_{m=1}^{N}\varphi(x^{-1}_n x_m)\bar{z_n}z_m\geq 0, \ \forall N\in\mathbb{N}, \forall x_1, \cdots, x_N\in G, \forall z_1, \cdots, z_N\in\mathbb{C}.
$$
\begin{theorem}\cite[Proposition 8.4.2]{Wolf}
	Let $\varphi$ be a positive-definite function. Then,
	\begin{enumerate}
		\item[1)] $\forall x\in G, |\varphi(x)|\leq\varphi(e)$.
		\item[2)] $\forall x\in G, \varphi(x^{-1})=\overline{\varphi(x)}$.
	\end{enumerate}
\end{theorem} 
If $\varphi$ is a spherical function, then $\varphi(e)=1$ \cite[Proposition 6.1.5]{Dijk2009}. 
  
Let $\mathcal{S}^+$ denote  the set of positive definite spherical functions for the Gelfand pair $(G,K)$. It carries a Plancherel measure which ensures the inversion formula for the spherical Fourier transform. Indeed, the spherical Fourier  transform extends to an isometric isomorphism $L^2(G//K)\longrightarrow L^2(\mathcal{S}^+)$ and 
$$f(x)=\int_{\mathcal{S}^+}\widehat{f}(\varphi)\varphi(x)d\varphi, \quad f\in L^2(G//K),\, x\in G.$$   

\subsection{The Stockwell transform on Gelfand pairs}
We recall from \cite{Dosseh2026} the definition of the Stockwell transform on a Gelfand pair together with the two properties,  the inversion formula and the Plancherel-type isometry, that will actually be used in the sequel.

Let $(G, K)$ be a Gelfand pair and $\alpha$ a topological automorphism of $G$. There exists a positive constant $\delta_\alpha$ such that $d(\alpha(x))=\delta_\alpha dx$. It follows,
$$\int_{G} f(\alpha(x))dx =\delta^{-1}_{\alpha}\int_{G} f(x)dx, \quad f\in \mathcal{C}_c(G//K).$$
\begin{definition}[{\cite[Definition 1]{Dosseh2026}}]
	Let $(G,K)$ be a Gelfand pair, let $\alpha$ be a topological automorphism of $G$ and let $\theta \in L^1 (G//K)\cap L^2 (G//K)$.
	The Stockwell transform of $f\in L^2(G//K)$ with respect to $\alpha$ and $\theta$ is
	\begin{equation}\label{ST-definition}
		S_{\theta,\alpha}f(t,\varphi)=\delta_\alpha^{\frac{1}{2}}\int_G f(x)\overline{\varphi (x)\theta (\alpha (t^{-1}x))}dx, \quad t\in G,  \varphi \in \mathcal{S}_b(G,K).
	\end{equation}
	\end{definition}
If we set $$\theta_{\alpha, \varphi, t}(x)=\varphi(x)\theta(\alpha(t^{-1}x)),$$ then formula  (\ref{ST-definition}) reads 
$$S_{\theta,\alpha}f(t,\varphi)=\delta_\alpha^{\frac{1}{2}}\langle f,\theta_{\alpha,\varphi,t}\rangle.$$

The space $\mathcal S^+$ is endowed with the Radon measure $d\varphi$ carrying the spherical Fourier inversion formula.

Set $$L^2(G\times\mathcal{S}^+)=\left\{f:G\times\mathcal{S}^+\to \mathbb{C} \text{ measurable such that } \int_{\mathcal{S}^+}\int_{G}|f(t, \varphi)|^2 dt d\varphi<\infty  \right\},$$ equipped with the norm
$$\|f\|_{L^2(G\times\mathcal{S}^+)}=\left(\int_{\mathcal{S}^+}\int_{G}|f(t, \varphi)|^2 dt d\varphi\right)^{\frac{1}{2}}.$$

\begin{theorem}[Inversion formula, {\cite[Equation (12)]{Dosseh2026}}]\label{thm:inversion}
	Let $\theta\in L^1(G//K)\cap L^2(G//K)$, $\theta\neq 0$. For every $f\in L^2(G//K)$,
	$$f(x)=\dfrac{\delta_\alpha^{\frac{1}{2}}}{\|\theta\|^2_{L^2}}\int_{\mathcal{S}^+}\int_{G} S_{\theta, \alpha}f(t, \varphi)\theta_{\alpha, \varphi, t}(x)\, dtd\varphi,\quad x\in G.$$
	In particular, when $\|\theta\|_{L^2}=1$,
	$$f(x)=\delta_\alpha^{\frac{1}{2}}\int_{\mathcal{S}^+}\int_{G} S_{\theta, \alpha}f(t, \varphi)\theta_{\alpha, \varphi, t}(x)\, dtd\varphi,\quad x\in G.$$
\end{theorem}

\begin{theorem}[Plancherel-type identity, {\cite[Theorem 7]{Dosseh2026}}]\label{thm:isometry}
	Let $(G,K)$ be a Gelfand pair and let $\theta,\eta, f, g \in L^2 (G//K)$. Then
	$$\langle S_{\theta,\alpha}f, S_{\eta,\alpha}g\rangle_{L^2(G\times\mathcal{S}^+)}=\langle f, g\rangle_{L^2(G//K)}\,\langle \eta,\theta\rangle_{L^2(G//K)}.$$
	In particular, if moreover $\|\theta\|_{L^2}=1$, then $S_{\theta,\alpha}$ is an isometry:
	\begin{eqnarray*}
		\|S_{\theta,\alpha}f\|_{L^2(G\times\mathcal{S}^+)}=\|f\|_{L^2}, \quad \forall f\in L^2(G//K).
	\end{eqnarray*}
\end{theorem}
\section{Main results}\label{main:results}
We state the main definition of this article. 
\begin{definition}\label{def:sobolev}
Let $\gamma: G\times\mathcal{S}^+\to \mathbb{R}_+$ be a positive measurable function and $s\in\mathbb{R}_+$. The Sobolev space associated with the Stockwell transform on the Gelfand pair $(G,K)$ is   the set
$$H_{\gamma, \theta,\alpha}^{s}(G//K)=\left\{ f\in L^2(G//K) : \int_{\mathcal{S}^+}\int_{G}(1+\gamma (t,\varphi)^2)^{s}|S_{\theta,\alpha}f(t, \varphi)|^2 dtd\varphi<\infty\right\}$$
endowed  with the norm $$ \|f\|_{H_{\gamma, \theta,\alpha}^{s}}=\left(\int_{\mathcal{S}^+}\int_{G}(1+\gamma (t,\varphi)^2)^{s}|S_{\theta,\alpha}f(t, \varphi)|^2 dtd\varphi \right)^{\frac{1}{2}}.$$
\end{definition}
 
We first show that the Sobolev-type space just introduced  carries a Hilbert space structure.
\begin{theorem}
	Let $(G,K)$ be a Gelfand pair and  let $\theta \in L^1 (G//K)\cap L^2(G//K)$ be such that $\|\theta\|_{L^2}=1$.
	The space $H_{\gamma, \theta,\alpha}^{s}(G//K)$  is a Hilbert space.
\end{theorem}
\begin{proof}
	Define $\Gamma: H_{\gamma, \theta,\alpha}^{s}(G//K)\to L^2(G\times\mathcal{S}^+)$ by $\Gamma f=(1+\gamma (\cdot,\cdot)^2)^{\frac{s}{2}}S_{\theta,\alpha}f(\cdot,\cdot)$. The linearity of $\Gamma$ is obvious. Furthermore,  we have
	\begin{align*}
		\|f\|_{H_{\gamma, \theta,\alpha}^{s}}^2&=\int_{\mathcal{S}^+}\int_{G}(1+\gamma (t,\varphi)^2)^{s}|S_{\theta,\alpha}f(t, \varphi)|^2 dtd\varphi\\
		&=\|(1+\gamma (\cdot,\cdot)^2)^{\frac{s}{2}}S_{\theta,\alpha}f(\cdot,\cdot)\|_{L^2(G\times\mathcal{S}^+)}^2\\
		&=\| \Gamma f\|_{L^2(G\times\mathcal{S}^+)}^2.
	\end{align*}
	So $\Gamma$ is a linear isometry, hence injective. In particular $\|\cdot\|_{H_{\gamma,\theta,\alpha}^s(G//K)}$ is the norm on $H_{\gamma,\theta,\alpha}^s(G//K)$, induced by the inner product
$$\langle f, h\rangle_{H^s_{\gamma,\theta,\alpha}}=\int_{\mathcal S^+}\int_{G}(1+\gamma(t,\varphi)^2)^s S_{\theta,\alpha}f(t,\varphi)\overline{S_{\theta,\alpha}h(t,\varphi)}\,dt\,d\varphi.$$

Let $(f_n)_n$ be a Cauchy sequence in $H_{\gamma,\theta,\alpha}^s(G//K)$. Since $s\geq 0$, we have $(1+\gamma(t,\varphi)^2)^s\geq 1$ for all $(t,\varphi)$, so by Theorem \ref{thm:isometry},
	$$\|f_n-f_m\|_{L^2}^2=\|S_{\theta,\alpha}(f_n-f_m)\|_{L^2(G\times\mathcal{S}^+)}^2\leq \|f_n-f_m\|_{H_{\gamma,\theta,\alpha}^s}^2.$$
Thus,  $(f_n)_n$ is also Cauchy in $L^2(G//K)$. Since $L^2(G//K)$ is complete, $(f_n)_n$ converges to some $f$ in $L^2(G//K)$.
	
On the other hand, $(\Gamma f_n)_n=\left((1+\gamma^2)^{\frac{s}{2}}S_{\theta,\alpha}f_n\right)_n$ is Cauchy in $L^2(G\times\mathcal{S}^+)$ as $\Gamma$ is an isometry, hence converges to some $g\in L^2(G\times\mathcal{S}^+)$. That is, $\left(S_{\theta,\alpha}f_n\right)_n$ converges to $(1+\gamma^2)^{-\frac{s}{2}}g$.
	
	Since $S_{\theta,\alpha}$ is an isometry $L^2(G//K)\to L^2(G\times\mathcal{S}^+)$, it is in particular continuous, so $S_{\theta,\alpha}f_n$ converges to $S_{\theta,\alpha}f$ in $L^2(G\times\mathcal{S}^+)$.  By uniqueness of the limit, $S_{\theta,\alpha}f=(1+\gamma^2)^{-\frac{s}{2}}g$ a.e., that is $g=\Gamma f$.
	
Since $\Gamma f=g\in L^2(G\times\mathcal S^+)$,  $f\in H_{\gamma,\theta,\alpha}^s(G//K)$ and
	$$\|f_n-f\|_{H_{\gamma,\theta,\alpha}^s}=\|\Gamma f_n-\Gamma f\|_{L^2(G\times\mathcal{S}^+)}=\|\Gamma f_n-g\|_{L^2(G\times\mathcal{S}^+)}\longrightarrow 0,$$
as $n$ goes to $\infty$.	Hence, every Cauchy sequence in $H_{\gamma,\theta,\alpha}^s(G//K)$ converges in $H_{\gamma,\theta,\alpha}^s(G//K)$, which is therefore a Hilbert space.
\end{proof}

Next, we show that increasing the smoothness parameter yields a decreasing scale of spaces, in the sense of continuous embeddings.
\begin{theorem} Let $(G,K)$ be a Gelfand pair and  let $\theta \in L^1 (G//K)\cap L^2(G//K)$ be such that $\|\theta\|_{L^2}=1$. If $\tau>\sigma\ge 0$,  then the  space $H_{\gamma, \theta,\alpha}^{\tau}(G//K)$ is continuously embedded into $H_{\gamma, \theta,\alpha}^{\sigma}(G//K)$ with $$\|f\|_{H_{\gamma, \theta,\alpha}^{\sigma}}\leqslant \|f\|_{H_{\gamma, \theta,\alpha}^{\tau}},\, \forall f\in H_{\gamma, \theta,\alpha}^{\tau}(G//K).$$
\end{theorem}
\begin{proof}
Since $(1+\gamma(t,\varphi)^2)\ge 1$, then for  $\tau>\sigma\ge 0$, we have  $(1+\gamma(t,\varphi)^2)^{\tau}\ge (1+\gamma(t,\varphi)^2)^{\sigma}$.
 The assertion follows. 
\end{proof}

We next establish an interpolation inequality, showing that the norm at an intermediate smoothness level is controlled by the norms at the two endpoint levels.
\begin{theorem}
Let $s_1, s_2$ be nonnegative real numbers satisfying $s_1<s_2$. Let $\zeta\in[0, 1]$, put $s=(1-\zeta)s_1+\zeta s_2$.  For all $f\in H_{\gamma, \theta,\alpha}^{s}(G//K)$, we have
	$$\|f\|_{H_{\gamma, \theta,\alpha}^{s}}\leq  \|f\|_{H_{\gamma, \theta,\alpha}^{s_1}}^{1-\zeta} \|f\|_{H_{\gamma, \theta,\alpha}^{s_2}}^{\zeta}.$$ 
\end{theorem}
\begin{proof}
	Let $\forall f\in H_{\gamma, \theta,\alpha}^{s}(G//K)$. We have
	\begin{align*}
		\|f\|_{H_{\gamma, \theta,\alpha}^{s}}^2&=\int_{\mathcal{S}^+} \int_{G}(1+\gamma (t,\varphi)^2)^{s}|S_{\theta,\alpha}f(t, \varphi)|^{2} dtd\varphi\\
		&=\int_{\mathcal{S}^+}\int_{G}(1+\gamma (t,\varphi)^2)^{[(1-\zeta)s_1+\zeta s_2]}|S_{\theta,\alpha}f(t, \varphi)|^{2\zeta+2(1-\zeta)} dtd\varphi\\
		&=\int_{\mathcal{S}^+}\int_{G}(1+\gamma (t,\varphi)^2)^{(1-\zeta)s_1}(1+\gamma (t,\varphi)^2)^{\zeta s_2}\times \\
		&|S_{\theta,\alpha}f(t, \varphi)|^{2\zeta} |S_{\theta,\alpha}f(t, \varphi)|^{2(1-\zeta)} dtd\varphi\\
		&\leq\left( \int_{\mathcal{S}^+}\int_{G}(1+\gamma (t,\varphi)^2)^{s_1}|S_{\theta,\alpha}f(t, \varphi)|^{2} dtd\varphi\right)^{1-\zeta}\times\\
		& \left(\int_{\mathcal{S}^+}\int_{G}(1+\gamma (t,\varphi)^2)^{s_2}|S_{\theta,\alpha}f(t, \varphi)|^{2} dtd\varphi \right)^{\zeta} \\
		&\leq \|f\|_{H_{\gamma, \theta,\alpha}^{s_1}}^{2(1-\zeta)} \|f\|_{H_{\gamma, \theta,\alpha}^{s_2}}^{2\zeta}.
	\end{align*}
 We deduce
 $$\|f\|_{H_{\gamma, \theta,\alpha}^{s}}\leq \|f\|_{H_{\gamma, \theta,\alpha}^{s_1}}^{1-\zeta}\|f\|_{H_{\gamma, \theta,\alpha}^{s_2}}^{\zeta}.$$
\end{proof}

In order to obtain embeddings into other Lebesgue spaces, we first establish the following auxiliary estimate, which bounds the Lebesgue norm of a function by a Lebesgue-type norm of its Stockwell transform.

\begin{lemma}\label{lemma:lpq}
Let $(G,K)$ be a Gelfand pair and let $\theta\in L^1(G//K)\cap L^2(G//K)$ with $\|\theta\|_{L^2}=1$.
Let $1\le p<\infty$ and $1\le q\le\infty$, with conjugate exponent $q'$, and set
$$N_{p,q}=\left(\int_{\mathcal{S}^+}\int_{G}\|\theta_{\alpha,\varphi,t}\|_{L^p(G)}^{q'}\,dt\,d\varphi\right)^{\frac{1}{q'}}$$
(an essential supremum over $(t,\varphi)$ if $q'=\infty$). If $N_{p,q}<\infty$, then for every $f\in L^2(G//K)$,
$$\|f\|_{L^{p}}\leq \delta_\alpha^{\frac{1}{2}}\,N_{p,q}\,\|S_{\theta,\alpha}f\|_{L^q(G\times\mathcal{S}^+)}.$$
\end{lemma}
\begin{proof}
Let $g\in\mathcal{C}_c(G//K)\subset L^2(G//K)$. Since $S_{\theta,\alpha}g(t,\varphi)=\delta_\alpha^{1/2}\langle g,\theta_{\alpha,\varphi,t}\rangle$, Theorem \ref{thm:isometry} (with $\eta=\theta$, $\|\theta\|_{L^2}=1$) gives 
$$\langle f,g\rangle=\langle S_{\theta,\alpha}f,S_{\theta,\alpha}g\rangle=\delta_\alpha^{\frac{1}{2}}\int_{\mathcal{S}^+}\int_G S_{\theta,\alpha}f(t,\varphi)\,\langle\theta_{\alpha,\varphi,t},g\rangle\,dt\,d\varphi .$$
By H\"older's inequality in $G$, $|\langle\theta_{\alpha,\varphi,t},g\rangle|\le\|\theta_{\alpha,\varphi,t}\|_{L^p}\|g\|_{L^{p'}}$, where $p'$ is the conjugate exponent of $p$. Then H\"older's inequality in $(t,\varphi)$ yields
$$|\langle f,g\rangle|\le \delta_\alpha^{\frac{1}{2}}\|g\|_{L^{p'}}\int_{\mathcal{S}^+}\int_G|S_{\theta,\alpha}f|\,\|\theta_{\alpha,\varphi,t}\|_{L^p}\,dt\,d\varphi\le\delta_\alpha^{\frac{1}{2}}N_{p,q}\|S_{\theta,\alpha}f\|_{L^q}\|g\|_{L^{p'}}.$$
Since $f$ is $K$-bi-invariant, we have 
$$\|f\|_{L^p}=\sup\{|\langle f,g\rangle| : g\in\mathcal{C}_c(G//K),\ \|g\|_{L^{p'}}\le 1\}$$ (valid for $1\le p<\infty$). Therefore, 
$$\|f\|_{L^{p}}\leq \delta_\alpha^{\frac{1}{2}}\,N_{p,q}\,\|S_{\theta,\alpha}f\|_{L^q(G\times\mathcal{S}^+)}.$$

\end{proof}

Building on the previous lemma, we now show that a suitable  condition on the weight yields a continuous embedding of the Sobolev-type space into a Lebesgue space.

\begin{theorem}\label{thm:Lp}
Let $(G,K)$ be a Gelfand pair and let $\theta\in L^1(G//K)\cap L^2(G//K)$ with $\|\theta\|_{L^2}=1$.
Let $s>0$, $\beta>2$ and set $p=\dfrac{2\beta}{\beta+2}\in(1,2)$, so that $\frac{1}{p}=\frac{1}{2}+\frac{1}{\beta}$. Assume that
$$(1+\gamma^2)^{-\frac{s}{2}}\in L^\beta(G\times\mathcal{S}^+)\quad\text{and}\quad N_p:=N_{p,p}<\infty$$
(with $N_{p,p}$ as in Lemma \ref{lemma:lpq}, i.e. with $q=p$). Then $H^s_{\gamma,\theta,\alpha}(G//K)$ is continuously embedded in $L^p(G//K)$ and
$$\|f\|_{L^p}\leq \delta_\alpha^{\frac{1}{2}}\,N_p\,\big\|(1+\gamma^2)^{-\frac{s}{2}}\big\|_{L^\beta(G\times\mathcal{S}^+)}\,\|f\|_{H^s_{\gamma,\theta,\alpha}},\quad f\in H^s_{\gamma,\theta,\alpha}(G//K).$$
\end{theorem}
\begin{proof}
Since $\frac{1}{p}=\frac{1}{2}+\frac{1}{\beta}$, the  H\"older inequality gives
\begin{align*}
\|S_{\theta,\alpha}f\|_{L^p(G\times\mathcal{S}^+)}&\le\|(1+\gamma^2)^{s/2}S_{\theta,\alpha}f\|_{L^2(G\times\mathcal{S}^+)}\,\|(1+\gamma^2)^{-s/2}\|_{L^\beta(G\times\mathcal{S}^+)}\\
&=\|f\|_{H^s_{\gamma,\theta,\alpha}}\,\|(1+\gamma^2)^{-\frac{s}{2}}\|_{L^\beta}.
\end{align*}
Lemma \ref{lemma:lpq} with $q=p$ concludes 
$\|f\|_{L^p}\le\delta_\alpha^{1/2}N_p\|S_{\theta,\alpha}f\|_{L^p}$.
\end{proof}

We now give a second embedding into Lebesgue spaces, based on the endpoint case $q=1$ of Lemma \ref{lemma:lpq}. It requires only that the window belongs to $L^p$ and that the inverse of the weight is square integrable, and it involves no condition on the family of analysing windows.
\begin{proposition}\label{prop:Lp-q1}
Let $\theta\in L^1(G//K)\cap L^2(G//K)\cap L^p(G//K)$ with $\|\theta\|_{L^2}=1$, $1\le p<\infty$, $s>0$, and assume $(1+\gamma^2)^{-s/2}\in L^2(G\times\mathcal{S}^+)$. Then $H^s_{\gamma,\theta,\alpha}(G//K)\hookrightarrow L^p(G//K)$ with
$$\|f\|_{L^p}\le\delta_\alpha^{\frac{1}{2}-\frac{1}{p}}\,\|\theta\|_{L^p}\,\big\|(1+\gamma^2)^{-\frac{s}{2}}\big\|_{L^2(G\times\mathcal{S}^+)}\,\|f\|_{H^s_{\gamma,\theta,\alpha}}.$$
\end{proposition}
\begin{proof}
Since $|\varphi|\le1$ and $\displaystyle\int_G|\theta(\alpha(t^{-1}x))|^p\,dx=\delta_\alpha^{-1}\|\theta\|_{L^p}^p$ (left invariance of the Haar measure and the change of variables $x\mapsto\alpha(x)$), we get $N_{p,1}\le\delta_\alpha^{-1/p}\|\theta\|_{L^p}$. Moreover, by Cauchy--Schwarz, $$\|S_{\theta,\alpha}f\|_{L^1}\le\|(1+\gamma^2)^{s/2}S_{\theta,\alpha}f\|_{L^2}\|(1+\gamma^2)^{-s/2}\|_{L^2}.$$
Finally, apply Lemma \ref{lemma:lpq} with $q=1$, $q'=\infty$ yields the result. 
\end{proof}

We finally turn to pointwise regularity. The following theorem states that, under a local domination condition on the window, the elements of the Sobolev-type space are continuous, and that  the space embeds continuously into the bounded continuous functions.

\begin{theorem}\label{thm:embedCb}
Let $(G,K)$ be a Gelfand pair with $G$ first countable, and let $\alpha$ be a topological automorphism of $G$ with $\alpha(K)=K$.
Let $\theta\in L^1(G//K)\cap L^2(G//K)$ be continuous with $\|\theta\|_{L^2}=1$, and let $s\ge0$. Consider the following hypotheses : 
\begin{itemize}
\item[(H)]  There exist a $K$-bi-invariant open neighbourhood $V$ of $x_0\in G$ and $h\in L^1(G\times\mathcal{S}^+)$ such that, for every $y\in V$ and for a.e. $(t,\varphi)$,
$$(1+\gamma(t,\varphi)^2)^{-s}\,|\theta(\alpha(t^{-1}y))|^2\le h(t,\varphi).$$
\item[(M)] $M:=\displaystyle\sup_{y\in G}\int_{\mathcal{S}^+}\int_G(1+\gamma(t,\varphi)^2)^{-s}|\theta(\alpha(t^{-1}y))|^2\,dt\,d\varphi<\infty.$
\end{itemize}
Then, 
\begin{enumerate}
\item[(a)] If (H) holds at $x_0$, then for every $f\in H^s_{\gamma,\theta,\alpha}(G//K)$ the integral
$$F(x)=\delta_\alpha^{\frac12}\int_{\mathcal{S}^+}\int_G S_{\theta,\alpha}f(t,\varphi)\,\varphi(x)\,\theta(\alpha(t^{-1}x))\,dt\,d\varphi$$
converges absolutely for $x\in V$, $F$ is continuous on $V$, and $F=f$ almost everywhere on $V$. In particular $f$ has a representative which is continuous at $x_0$.
\item[(b)] If (H) holds at every point of $G$ and (M) holds, then every $f\in H^s_{\gamma,\theta,\alpha}(G//K)$ has a representative $F\in\mathcal{C}_b(G//K)$, and
$$\|F\|_\infty\le\delta_\alpha^{\frac12}M^{\frac12}\|f\|_{H^s_{\gamma,\theta,\alpha}}.$$
Thus $H^s_{\gamma,\theta,\alpha}(G//K)\hookrightarrow\mathcal{C}_b(G//K)$ continuously.
\end{enumerate}
\end{theorem}

\begin{proof}
Let $f\in H^s_{\gamma,\theta,\alpha}(G//K)$. Set $\theta_x(t,\varphi)=\theta(\alpha(t^{-1}x))$.
\begin{enumerate}
\item[\textnormal{(a)}] For $x\in V$, using $|\varphi(x)|\le1$, the Cauchy-Schwarz inequality and (H), we have
\begin{align*}
\int_{\mathcal{S}^+}\int_G|S_{\theta,\alpha}f(t,\varphi)||\theta_x(t,\varphi)|\,dt\,d\varphi&\le\|(1+\gamma^2)^{s/2}S_{\theta,\alpha}f\|_{L^2}\left(\int_{\mathcal{S}^+}\int_G (1+\gamma(t,\varphi)^2)^{-s}|\theta_x(t,\varphi)|^2\,dt\,d\varphi\right)^{\frac12}\\
&(\text{using } |\varphi(x)|\le1 \text{ and the Cauchy-Schwarz inequality})\\
&\le\|f\|_{H^s_{\gamma,\theta,\alpha}}\|h\|_{L^1}^{\frac{1}{2}}<\infty\\
& (\text{using (H)}).
\end{align*}

 Let $x_0\in V$ and let  $(x_n)$ be a sequence which converges to $x_0$. Since $V$ is open, $x_n\in V$ for $n$ large. Put $I_n(t,\varphi)=S_{\theta,\alpha}f(t,\varphi)\varphi(x_n)\theta(\alpha(t^{-1}x_n))$. Since there are countably many $x_n$, for a.e. $(t,\varphi)$ and all $n$, by (H), we have 
\begin{align*}
|I_n(t,\varphi)|&\le|S_{\theta,\alpha}f(t,\varphi)|\,|\theta_{x_n}(t,\varphi)|\\
&\le (1+\gamma^2(t,\varphi))^{s/2}|S_{\theta,\alpha}f(t,\varphi)|\,h(t,\varphi)^{\frac{1}{2}}=:\mathcal{G}(t,\varphi).
\end{align*}
Let us now set 
$$\mathcal{G}(t,\varphi)=(1+\gamma^2(t,\varphi))^{s/2}|S_{\theta,\alpha}f(t,\varphi)|\,h(t,\varphi)^{\frac{1}{2}}.$$
We have $\mathcal{G}\in L^1(G\times\mathcal{S}^+)$ because $(1+\gamma^2)^{s/2}S_{\theta,\alpha}f\in L^2(G\times\mathcal{S}^+)$ and $h^{1/2}\in L^2(G\times\mathcal{S}^+)$.  Moreover, for every $(t,\varphi)$, $I_n(t,\varphi)\to I_0(t,\varphi)$, because $\varphi$ and $\theta$ are continuous and so are $\alpha$ and the group operations. The Lebesgue dominated convergence theorem implies that $F(x_n)$ converges to $F(x_0)$. Since   $G$ is first countable, $F$ is continuous at $x_0$.

Now let us show that $F=f$ a.e. on $V$. First, $F$ is $K$-bi-invariant on $V$: left $K$-invariance follows from the change of variables $t\mapsto kt$ together with $S_{\theta,\alpha}f(kt,\varphi)=S_{\theta,\alpha}f(t,\varphi)$ (because $f$ and $\varphi$ are $K$-bi-invariant); right $K$-invariance follows from $\varphi(xk)=\varphi(x)$, $\alpha(K)=K$ and the right $K$-invariance of $\theta$. Now let $g\in\mathcal{C}_c(G//K)$ with $\operatorname{supp}g\subset V$. 
We have
\begin{align*}
\int_G F(x)\,\overline{ g(x)}\,dx&=\delta_\alpha^{\frac12}\int_G \int_{\mathcal{S}^+}\int_G S_{\theta,\alpha}f(t,\varphi)\,\varphi(x)\,\theta(\alpha(t^{-1}x))\,dt\,d\varphi\overline{ g(x)}\,dx\\
&= \delta_\alpha^{\frac12}\int_{\mathcal{S}^+}\int_G S_{\theta,\alpha}f(t,\varphi)\,\langle\theta_{\alpha,\varphi,t},g\rangle\,dt\,d\varphi\\
& \text{(by Fubini's theorem)}\\
&=\langle f,g\rangle.
\end{align*}

Hence $\displaystyle\int_V(f-F)\overline g=0$ for all such $g$. Since  $f-F$ is $K$-bi-invariant and locally integrable on $V$, $f=F$ a.e. on $V$. This proves (a).

\item[\textnormal{(b)}] Under (H) everywhere, $F$ is defined and continuous on $G$ and equals $f$ a.e. Under (M), for every $x\in G$, the inequality 
\begin{align*}
\int_{\mathcal{S}^+}\int_G|S_{\theta,\alpha}f(t,\varphi)||\theta_x(t,\varphi)|\,dt\,d\varphi&\le\|(1+\gamma^2)^{s/2}S_{\theta,\alpha}f\|_{L^2}\left(\int_{\mathcal{S}^+}\int_G (1+\gamma(t,\varphi)^2)^{-s}|\theta_x(t,\varphi)|^2\,dt\,d\varphi\right)^{\frac12}
\end{align*}
implies
$$|F(x)|\le\delta_\alpha^{\frac12}\int_{\mathcal{S}^+}\int_G|S_{\theta,\alpha}f(t,\varphi)||\theta_x(t,\varphi)|\,dt\,d\varphi\le\delta_\alpha^{\frac12}M^{\frac12}\|f\|_{H^s_{\gamma,\theta,\alpha}}.$$
Therefore, 
$$\|F\|_\infty\le\delta_\alpha^{\frac12}M^{\frac12}\|f\|_{H^s_{\gamma,\theta,\alpha}},$$
 and the  following continuous embedding follows : $H^s_{\gamma,\theta,\alpha}(G//K)\hookrightarrow\mathcal{C}_b(G//K)$.

\end{enumerate}
\end{proof}

To illustrate the theory, we now work out the construction of the Sobolev-type space explicitly on the Euclidean motion group.

\begin{example}[The Euclidean motion group]\label{ex:En}
{\rm
Let $n\ge2$ and $(G,K)=(E(n),SO(n))$, where $E(n)=\mathbb{R}^n\rtimes SO(n)$ with product and inverse
$$(y,k)(y',k')=(y+k\cdot y',\,kk'),\quad (y,k)^{-1}=(-k^{-1}y,\,k^{-1}).$$
The Haar measure of $G$ is $dy\,dk$, where $dy$ is the Lebesgue measure on $\mathbb{R}^n$ and $dk$ is the normalized Haar measure of $K$.  A function on $G$ is $K$-bi-invariant if and only if it is of the form $(y,k)\mapsto f(|y|)$, i.e.\ it is a radial function on $\mathbb{R}^n$ \cite[Subsection 7.1]{Dijk2009}. We identify $L^p(G//K)$ with the radial functions of $L^p(\mathbb{R}^n)$.

\smallskip
The bounded spherical functions are indexed by $\lambda\ge0$ :
$$\varphi_\lambda(y)=\Gamma\!\left(\tfrac{n}{2}\right)2^{\frac{n-2}{2}}\,\frac{J_{\frac{n-2}{2}}(\lambda|y|)}{(\lambda|y|)^{\frac{n-2}{2}}},\quad |\varphi_\lambda|\le1,\quad\varphi_\lambda(0)=1,$$
where $J_\nu$ is the Bessel function of the first kind. They are real-valued, $\mathcal{S}^+\cong[0,\infty)$, and the spherical Fourier transform $\widehat f(\lambda)=\displaystyle\int_{\mathbb{R}^n}f(y)\varphi_\lambda(y)\,dy$ is the Hankel transform. With the Haar measure $dy$, the Plancherel measure is $d\varphi_\lambda=c_n\lambda^{n-1}d\lambda$, where $c_n>0$ depends only on $n$ and on
  the normalization of the Haar measure. 
   
\smallskip
 For $a>0$ let $\alpha_a(y,k)=(ay,k)$. Then $\alpha_a$ is a topological automorphism of $E(n)$ with $\alpha_a(K)=K$ (it fixes $K$ pointwise), and $d(\alpha_a(x))=a^n dx$, i.e.\ $\delta_{\alpha_a}=a^n$. Let
$$\theta(y,k)=\theta(y)=\pi^{-\frac n4}e^{-\frac{|y|^2}{2}},$$
which is radial, continuous, belongs to $L^1\cap L^2$, and satisfies $\|\theta\|_{L^2}=1$. For $t=(u,h)$ and $x=(y,k)$ in $G$,
$$\alpha_a(t^{-1}x)=\big(a\,h^{-1}(y-u),\,h^{-1}k\big),\quad\text{hence}\quad\theta(\alpha_a(t^{-1}x))=\pi^{-\frac n4}e^{-\frac{a^2|y-u|^2}{2}},$$
which depends neither on $h$ nor on $k$. Since $\varphi_\lambda$ and $\theta$ are real and the integral over $K$ is trivial, for a radial $f$ we get
\begin{equation}\label{eq:ST-Enk}
S_{\theta,\alpha_a}f(u,h,\lambda)=S_{\theta,\alpha_a}f(u,\lambda)=a^{\frac n2}\pi^{-\frac n4}\int_{\mathbb{R}^n}f(y)\,\varphi_\lambda(y)\,e^{-\frac{a^2|y-u|^2}{2}}\,dy,\quad u\in\mathbb{R}^n,\ \lambda\ge0,
\end{equation}
a windowed Hankel-type transform.

\smallskip
 Take the weight $\gamma(t,\lambda)=\lambda$, which depends only on the spherical variable. Definition \ref{def:sobolev} gives
$$H^s_{\gamma,\theta,\alpha_a}\big(E(n)/\!/SO(n)\big)=\Big\{f\text{ radial on }\mathbb{R}^n:c_n\int_0^\infty\int_{\mathbb{R}^n}(1+\lambda^2)^s\,|S_{\theta,\alpha_a}f(u,\lambda)|^2\,du\,\lambda^{n-1}d\lambda<\infty\Big\}.$$

\smallskip
We check the hypotheses of Theorem \ref{thm:embedCb}. Fix $\rho\ge0$, $r>0$ and the $K$-bi-invariant open set $V=\{(y,k):\ ||y|-\rho|<r\}$; every point of $G$ has such a neighbourhood. For $y\in V$ and $u\in\mathbb{R}^n$ we have $|y-u|\ge||u|-\rho|-r$ whenever $||u|-\rho|\ge r$, hence, with $t=(u,h)$,
$$|\theta(\alpha_a(t^{-1}y))|^2=\pi^{-\frac{n}{2}}e^{-a^2|y-u|^2}\le\pi^{-\frac{n}{2}}\kappa(u),\quad
\kappa(u)=\begin{cases}1,&||u|-\rho|<r,\\ e^{-a^2(||u|-\rho|-r)^2},&||u|-\rho|\ge r.\end{cases}$$
The function $\kappa$ is integrable on $\mathbb{R}^n$ (a bounded annulus and a Gaussian tail). Put $h(t,\lambda)=\pi^{-n/2}\kappa(u)(1+\lambda^2)^{-s}$: it satisfies (H) on $V$, and it is integrable on $G\times\mathcal{S}^+$ as soon as $\displaystyle\int_0^\infty(1+\lambda^2)^{-s}\lambda^{n-1}d\lambda<\infty$. Moreover, for every $y\in G$, the change of variables $w=t^{-1}y$ and the unimodularity of $G$ give $\int_G|\theta(\alpha_a(t^{-1}y))|^2dt=\delta_{\alpha_a}^{-1}=a^{-n}$, so that
$$M=a^{-n}c_n\int_0^\infty(1+\lambda^2)^{-s}\lambda^{n-1}\,d\lambda ,$$
which is finite if and only if $s>n/2$. Theorem \ref{thm:embedCb} therefore shows that, for $s>n/2$, every element of $H^s_{\gamma,\theta,\alpha_a}(E(n)/\!/SO(n))$ has a representative $F$ which is a bounded continuous radial function on $\mathbb{R}^n$, with
$$\|F\|_\infty\le\Big(c_n\int_0^\infty(1+\lambda^2)^{-s}\lambda^{n-1}d\lambda\Big)^{\frac{1}{2}}\|f\|_{H^s_{\gamma,\theta,\alpha_a}}.$$
 The condition $s>n/2$ is the classical threshold for the embedding of $H^s(\mathbb{R}^n)$ into the bounded continuous functions \cite{SteinWeiss} where
  $$
H^s(\mathbb{R}^n) = \left\{ u \in \mathcal{S}'(\mathbb{R}^n) : (1 + |\xi|^2)^{s/2} \widehat{u}(\xi) \in L^2(\mathbb{R}^n) \right\},
$$ 
$\mathcal{S}'(\mathbb{R}^n)$ being the space of 
 tempered distributions  and   $\widehat{u}$ is the Fourier transform of $u$.

}
\end{example}

 \section{Conclusion}

In this paper we have defined Sobolev-type spaces on a Gelfand pair through the Stockwell transform, with the smoothness measured by a weight acting on the time-frequency variables rather than on the frequency alone. The Plancherel-type identity for the Stockwell transform makes these spaces Hilbert spaces, and the usual structural properties of a Sobolev scale hold: continuous embeddings between different smoothness levels and an interpolation inequality. We also obtained embeddings into Lebesgue spaces and into the bounded continuous functions. The embedding into bounded continuous functions rests on a domination condition on the window and requires no equicontinuity of the family of spherical functions. On the Euclidean motion group, the construction leads to a windowed Hankel transform and to the classical smoothness threshold for continuity.

Finally, the same construction could be attempted for the Stockwell transform in related frameworks, such as hypergroups or commutative triples, and applied to the regularity of localization operators.

\end{document}